\documentclass{article}
\usepackage{chicago}
\usepackage{epsfig}
\usepackage{amsmath,amssymb}
\title{The existence of maximum likelihood estimates in the
Bradley-Terry model and its extensions}
\author{Kenneth Butler and John T.\ Whelan\thanks{jtw24@cornell.edu}\\
Dalhousie University and University of Bern}

\newcommand{\tre}{\trianglerighteq}
\newcommand{\D}{\displaystyle}

\renewcommand{\baselinestretch}{1.8}

\begin{document}
\maketitle

\begin{abstract}
In the Bradley-Terry model for paired comparisons, and its extensions
to include order effects and ties, the maximum likelihood estimates
of probabilities of certain outcomes can be 0 or 1 under certain data
configurations. This poses
problems for standard estimation methods. In this paper, we give
algorithms for identifying the outcomes with estimated probability 0
or 1, and indicate how the remaining probabilities may be estimated
and summarized.
\end{abstract}

\section{Introduction} \label{intro}

A paired-comparison experiment is used to assess the relative worths
of $t$ objects when they can only be compared two at a time, and when
the result of such a comparison yields only the result that one of the
objects was preferred to the other (or, possibly, that the two objects
in the comparison were equally preferable). The most commonly-used
model for data from such an experiment is that of \citeN{bradterr52},
with extensions by \citeN{davidson70} and \citeN{davibeav77}. These
models yield estimates of the probability that one object will be
preferred to another.  \citeN{davifarq76} give an extensive
bibliography of Bradley-Terry and related models.

Applications of the Bradley-Terry model include taste-testing, in
which each comparison is carried out by a judge, and sporting
competition, in which the objects are players or teams and comparisons
represent games between them, with preferences corresponding to wins
and losses (or ties).  In this paper, we use for definiteness the
terms ``teams'' and ``games'' to refer to the objects and comparisons,
though our results apply more widely.

A drawback of the Bradley-Terry model is that it is possible for
estimated probabilities to equal 0 or 1, and therefore the maximum
likelihood estimates of the parameters are not guaranteed to exist (in
the sense of being in the interior of the parameter space).  This was
noted by \citeN{ford57}, who gave a necessary and sufficient condition
for the estimated probabilities to lie strictly between 0 and 1 in
Bradley and Terry's original model.

The Bradley-Terry model and its extensions can be viewed as special
cases of (binomial and trinomial) logistic regression models.
\citeN{albeande84} and \citeN{santduff86} give conditions which
determine whether the maximum likelihood estimates in logistic
regression are finite or infinite, and give linear programs and other
methods for determining the effect of any nonexistent parameter
estimates on the estimated probabilities.

In this paper, we use the results of \citeANP{albeande84} and
\citeANP{santduff86} to show that a simple procedure, based on ideas
of graph theory, leads to a complete analysis of a paired comparison
experiment using the original Bradley-Terry model and one of its
extensions, even when the maximum likelihood estimates do not all
exist. We also demonstrate that variations on this procedure apply to
the other extensions of the model.  In Section~\ref{secbt}, we review
the Bradley-Terry model and its extensions; in Section~\ref{secnotie}
we consider estimation when the only possible outcomes are wins and
losses, and in Section~\ref{sectie} we consider what happens when ties
are included in the model. Finally, we discuss the presentation of
results from our methods in Section~\ref{secpres} and give some
illustrative examples in Section~\ref{secex}.

\section{Bradley-Terry model and extensions}\label{secbt}

\subsection{Bradley and Terry's original model}\label{secbtone}

The \citeN{bradterr52} model assumes that each team has a ``worth'' or
``strength'' $\pi_i \ge 0$ reflecting the $i$-th team's tendency to
win (or lose). The $\pi_i$ are determined only up to ratios, and so in
practice an additional condition, such as $\sum_{i=1}^t \pi_i=1$ or
$\pi_t=1$, is imposed to force the maximum likelihood estimates to be
well-defined.  The Bradley-Terry model asserts that the probability
$p_{ij}$ of team $i$ defeating team $j$ is
\begin{equation}
p_{ij}=\frac{\pi_i}{\pi_i+\pi_j};
\label{btorig}
\end{equation}
equivalently, $\pi_i/\pi_j$ represents the odds in favour of $i$
defeating $j$. (The Bradley-Terry model is therefore a
``proportional-odds'' model.)

If each pair of teams $i$ and $j$ play $n_{ij}$ games against each
other with team $i$ winning $y_{ij}$ of them, and all games are
assumed independent, the likelihood is seen to be

$$L = \prod_{i<j} \frac{\pi_i^{y_{ij}} \pi_j^{n_{ij}-y_{ij}}}
{(\pi_i+\pi_j)^{n_{ij}}}.$$
\citeN{fienberg79} shows that the
likelihood equations require the observed and expected numbers of wins
for each team to be equal (where the expected number of wins is the
sum of win probabilities over games played by the team).

\subsection{An order effect}\label{secbttwo}
In taste testing, an object may be more likely to be preferred simply
because it is presented first rather than second; in sports, a team is
typically more likely to defeat another when playing at home. This
``order effect'' or ``home field advantage'' can be included in the
model by assuming the probability of $i$ defeating $j$ in that order
(that is, with team $i$ at home) to be
\begin{equation}
 p_{ij}= \frac{ \gamma \pi_i }{ \gamma\pi_i + \pi_j }. 
\label{eqhomeprob}
\end{equation}

The additional parameter $\gamma \ge 0$ represents the advantage of
being first (playing at home); it is assumed to be the same for all
teams. This model was proposed by \citeN{davibeav77}. As in
Section~\ref{secbtone}, the likelihood equations equate the observed
and expected wins for each team; the additional equation for $\gamma$
equates the observed and expected wins for the home team over all
games played.

\subsection{Team-specific order effects}\label{secbtthree}
It is often the case that a single home field (order) effect is
consistent with the data, though this may be because there is
insufficient data to prove that a separate home field effect is
necessary for each team. Nonetheless, it is straightforward to
generalize (\ref{eqhomeprob}) to allow team-specific home field
effects. This is most easily done by observing that we are effectively
treating ``team $i$ at home'' and ``team $i$ on the road'' as
completely separate entities. Letting $\pi_{iH}$ denote the strength
of team $i$ at home and $\pi_{iV}$ the strength of team $i$ on the
road (as ``visiting team''), the probability $p_{ij}$ of team $i$
defeating team $j$ with $i$ at home can be modelled as

\begin{equation}
 p_{ij} = \frac{\pi_{iH} }{ \pi_{iH}+\pi_{jV} }
\label{allhomeprob}
\end{equation}
and the likelihood equations now equate observed and expected home and
road wins for each team (a total of $2t$ equations for the $2t$
parameters, though one equation and one parameter are redundant).

\subsection{Ties and the Bradley-Terry model}\label{secbtfour}
If one of the possible outcomes of a comparison is ``no preference''
(that is, ``tie''), then this should be included in the model. One
straightforward approach is simply to count a tie as half a win and
half a loss and use (\ref{btorig}), (\ref{eqhomeprob}) or
(\ref{allhomeprob}) as appropriate. This violates the binomial
assumption underlying those models, but can work well in practice if
the goal is simply to estimate the strengths of the teams (rather than
estimating probabilities).  \citeN{davidson70} proposed a
generalization of the Bradley-Terry model that allows for a tie as a
third outcome. Letting $p_{ijk}$ denote the probability of outcome
type $k$ ($k=1$ if $i$ wins, $k=2$ if $j$ wins, $k=0$ if a tie), the
generalization is
\begin{equation}
\begin{array}{rcl}
p_{ij1} &=& \frac{\D \pi_i }
  {\D \pi_i + \pi_j + \nu \sqrt{\pi_i \pi_j}}\\
p_{ij2} &=& \frac{\D \pi_j }
  {\D \pi_i + \pi_j + \nu \sqrt{\pi_i \pi_j}}\\
p_{ij0} &=& \frac{\D \nu \sqrt{\pi_i \pi_j} }
{\D \pi_i + \pi_j + \nu \sqrt{\pi_i \pi_j}}
\end{array}
\label{eqtieprob}
\end{equation}

The additional parameter $\nu\ge 0$ reflects the tendency for ties to
happen, and is assumed to be the same for all teams.  This model has
two desirable properties: the conditional probability of $i$ defeating
$j$, given a non-tie, is as in (\ref{btorig}), and the probability of
a tie, for fixed $\nu$, is maximum when $\pi_i=\pi_j$. The likelihood
equation for $\nu$ equates the observed and expected ties over all
games; the remaining equations now equate the observed and expected
{\em points} for each team, where one point is given for a win and
half a point for a tie (or equivalent such as two points for a win and
one for a tie).

\citeN{davibeav77} show how ties and order effects may be modelled
simultaneously; it is also possible to combine the above model with
(\ref{allhomeprob}).

\subsection{Team-specific tie effects}\label{secbtfive}
Model (\ref{eqtieprob}) assumes that all teams have the same tendency
to tie games, and this may not always be reasonable. For example,
\citeN{joe90} found that some chess players tend to draw more matches
than others, even after allowing for the strength of opposition faced.
In this case, one can replace $\nu$ by a $\nu_i$ for each team, as
follows:
\begin{equation}
\begin{array}{rcl}
p_{ij1} &=& \frac{\D \pi_i }
  {\D \pi_i + \pi_j + \sqrt{\nu_i \nu_j \pi_i \pi_j}}\\
p_{ij2} &=& \frac{\D \pi_j }
  {\D \pi_i + \pi_j + \sqrt{\nu_i \nu_j \pi_i \pi_j}}\\
p_{ij0} &=& \frac{\D \sqrt{\nu_i \nu_j \pi_i \pi_j}}
  {\D \pi_i + \pi_j + \sqrt{\nu_i \nu_j \pi_i \pi_j}}
\end{array}
\label{alltieprob}
\end{equation}
with the likelihood equations now equating the observed and expected
points and ties (and therefore wins) for each team.  This model can
also be combined with one containing order effects.

\subsection{Numerical methods}
\citeN{ford57} describes the following procedure for the model
(\ref{btorig}): the likelihood equations set, for team $i$,
$$ \sum_{j=1}^t \frac{n_{ij} \pi_i }{\pi_i+\pi_j } = \sum_{j=1}^t
y_{ij},$$
and therefore
$$ \pi_i= \sum_{j=1}^t y_{ij} \Big/ \sum_{j=1}^t \frac{ n_{ij} }
{ \pi_i+\pi_j }.$$
This can be treated as a fixed-point iteration, using the current
parameter values on the right-hand side to produce an updated value
for the parameter on the left.  Ford shows that, if the maximum
likelihood estimates $\hat{\pi}_i$ all exist, the values from the
iterative procedure will converge to the maximum likelihood estimates.
Convergence is generally slow, but the procedure is reliable and easy
to program. The same idea can be used in the extensions to the
Bradley-Terry model.

Several other iterative procedures are available, such as those
described and compared in \citeN[Section 3.3]{butler97}. These also
assume that the parameter estimates $\hat\pi_i$ all exist.  They tend
to require fewer iterations than Ford's method, but have greater
computational cost per iteration so that an overall improvement in
speed is not guaranteed.

However, the existence of the $\hat\pi_i$ 
is not normally apparent,
and if one or more of the $\hat\pi_i$ is actually infinite or
zero, an iterative procedure will only converge slowly, if at all.
Indeed, the usual reason for such a procedure to reach its iteration
limit is that some parameter estimates do not exist.

\section{Estimation with win/loss outcomes}\label{secnotie}

\subsection{Introduction}
When parameter estimates may not exist, a much more satisfactory
approach is first to identify the nonexistent ones (and hence the
fitted probabilities that are 1 or 0), and only then to apply an
iterative procedure if necessary. \citeN{ford57} showed that, in
(\ref{btorig}), the estimates will all exist if and only if there is
no division of the teams into two groups A and B such that every game
between a team in A and a team in B was won by the team in A. More
generally, in logistic regression, some observations may have fitted
probabilities that are 0 or 1, and it is desired to isolate and remove
these observations before applying an iterative procedure to estimate
the remainder of the probabilities; this is the problem addressed by
\citeN{albeande84} and \citeN{santduff86}.

\subsection{Simple Bradley-Terry model}\label{secestone}

The model (\ref{btorig}) is written as a logistic regression model by
letting $\beta_l=\log \pi_l$; \citeN{albeande84} showed that the
existence of maximum likelihood estimates in such a model was
connected with the existence of certain solutions to the corresponding
logistic discrimination problem.  In the present context, that problem
is defined by requiring
\begin{equation}
  \beta_{i}-\beta_{j} \ge 0
\label{eqsep}
\end{equation}
if there exists at least one game in which $i$ has defeated $j$.

\citeANP{albeande84} and \citeANP{santduff86} then characterize the
cases in which maximum likelihood estimates fail to exist by defining
``complete separation'', ``quasi-complete separation'' and
``overlap''. For complete separation, there must exist $\beta_i$
satisfying (\ref{eqsep}) with strict inequality for all observations.
In this case, all the fitted probabilities (for outcomes which
actually occurred) and the maximized likelihood itself are 1. For
quasi-complete separation, there must exist $\beta_i$ satisfying
(\ref{eqsep}) with at least one strict inequality and at least one
equality. In this case, the maximum of the likelihood is less than 1,
but some of the fitted probabilities are 1. Overlap occurs otherwise
(that is, when the only solution to (\ref{eqsep}) is with equality and
all the $\beta_i$ are equal); in this case, it is shown that all the
maximum likelihood estimates all exist and all the estimated
probabilities lie strictly between 0 and 1. (Note that the $\beta_i$
in a sense play a dual r\^{o}le: for determining the separation
properties of the data, we are concerned only with the existence of
non-equal $\beta_i$ solving the discrimination problem; the $\beta_i$
in the likelihood itself for our regression problem represent the
relative strengths of the teams, and the values will be different in
the two contexts.  For clarity, we will, as much as possible, work
with $\beta$ when discussing the question of separation and $\pi$ when
dealing with maximum likelihood estimates.)

For any two teams $i$ and $j$, define the relation $\ge$ by $i \ge j$
if there is a game in which $i$ has defeated $j$. (Note that we can
have $i \ge j$ and $j \ge i$, if for instance $i$ and $j$ have played
twice and won one game each.) For consistency in the sequel, we also
define $i \ge i$ to be true for each team $i$.

Define the relation $\trianglerighteq$ by $i \trianglerighteq j$ if $i
\ge j$ or if there exist teams $k_1, k_2, \ldots k_r$ such that $i \ge
k_1 \ge \cdots \ge k_r \ge j$.  Now consider the $t$ teams as vertices
of a directed graph $G$ which has an edge from vertex $i$ to vertex
$j$ if $i \ge j$. Determining whether the relation $\trianglerighteq$
holds for a pair of teams $k$ and $l$ is equivalent to determining
whether there is a path in the directed graph from the vertex
representing $k$ to that representing $l$ -- that is, whether there is
an edge from $k$ to $l$ in the {\em transitive closure} $\bar{G}$ of
$G$.  The transitive closure of a graph can be computed in $O(t^3)$
time by an algorithm due to Floyd and Warshall -- see, for example,
\citeN[p.\ 132]{papastei98}.

For each pair of teams $i$ and $j$, $i \trianglerighteq j$ or not and
$j \trianglerighteq i$ or not. Define the relations
$\cong,\gg,\gtrless$ as follows:
\begin{eqnarray*}
i \cong j & \mbox{if } 
& i \trianglerighteq j \mbox{ and } j \trianglerighteq i\\
i \gg j & \mbox{if } & i \trianglerighteq j \mbox{ only}\\
i \gtrless j & \mbox{if } 
& \mbox{neither } i \trianglerighteq j 
  \mbox{ nor } j \trianglerighteq i
\end{eqnarray*}
It is possible, though redundant, to define $i \ll j$ if and only if
$j \gg i$. Further, $i \cong i$ since $i \ge i$ and hence $i
\trianglerighteq i$.

It is easily seen that $\cong$ is an equivalence relation and that
$\gg$ is an ordering. We assert that when $i \cong j$, the fitted
probability of $i$ defeating $j$ lies strictly between 0 and 1, when
$i \gg j$, the fitted probability is 1, and when $i \gtrless j$, the
fitted probability is arbitrary since there is no basis for comparison
between $i$ and $j$.

To apply the Albert-Anderson-Santner-Duffy theory, we assign $\beta_i$
such that $\beta_i=\beta_j$ if $i \cong j$. For $i \gg j$, we assign
values such that $\beta_i > \beta_j$.  These assignments are possible
since $\cong$ is an equivalence relation and $\gg$ is an ordering.
Provided that all teams have played at least one game, this permits
values to be assigned to all $\beta_i$. In each case, the numerical
values of the $\beta_i$ are arbitrary; it is only the ordering that is
of importance.

To show that this correctly identifies complete separation,
quasi-complete separation, and overlap, note first that $i \cong j$ if
and only if $i \trianglerighteq j$ and $j \trianglerighteq i$. Thus
there exist (possibly empty) strings of teams such that $i \ge k_1 \ge
\cdots \ge k_r \ge j$ and $j \ge l_1 \ge \cdots \ge l_s \ge i$. The
discrimination equation (\ref{eqsep}) implies that $\beta_i \ge
\beta_{k_1} \ge \cdots \ge \beta_j$ and $\beta_j \ge \beta_{l_1} \ge
\cdots \ge \beta_i$. This can only be satisfied by $\beta_i=\beta_j$.
Thus $\beta_i=\beta_j$ if and only if $i \cong j$. When $i \gg j$,
however, $\beta_i>\beta_j$; at least one of the inequalities in the
chain must be strict since otherwise $j \trianglerighteq i$ as well.
Finally, if $i\gtrless j$, neither $i \cong j$ nor $j \cong i$, which
means that no chain of inequalities between $\beta_i$ and $\beta_j$ in
either direction can be inferred from the discrimination equation.
Thus there exist solutions with either sign of $\beta_i-\beta_j$.
(This is only possible if $i$ and $j$ did not play one another, and
thus is only relevant to the probabilities of outcomes of games which
did not occur.)  It follows that if no pair of teams is related by
$\cong$, the $\beta_i$ may all be chosen to be different, and we have
complete separation. If there exist teams $i$ and $j$ with $i \cong j$
and also teams $k$ and $l$ with $k \gg l$ (or $k \gtrless l$), we take
$\beta_i=\beta_j$ and $\beta_{k} >\beta_{l}$, indicating
quasi-complete separation.  Finally, if all teams are related by
$\cong$, then $\beta_i=\beta_j$ for all $i,j$, and so we have overlap.

The identification of outcomes with 0 or 1 probabilities allows us to
replace the probability of a win for $i$ in a game against $j$ in
(\ref{btorig}) with
$$ p_{ij} = \left\{
\begin{array}{rcl}
\displaystyle{\frac{\pi_i}{\pi_i+\pi_j}} & \mbox{if } & i \cong j\\
1 & \mbox{if } & i \gg j\\
0 & \mbox{if } & j \gg i\\
\mbox{ arbitrary} & \mbox{if } & i \gtrless j.
\end{array}
\right.
$$
An ``arbitrary'' probability is one whose value is not determined
by the data (i.e., there exist maximum likelihood solutions with all
values, between and including 0 and 1, for the probability in
question).  The effect of the modified definition of probability on
the maximum likelihood equations is to remove from the calculation any
games whose outcomes have maximum likelihood probabilities of 0 or 1.
The resulting strength estimates $\hat\pi_i$ are all guaranteed to
exist and can be found by methods such as those discussed in
Section~\ref{secbtfive}.  Maximum likelihood estimates of unit (or
vanishing) probability are now associated with a $\gg$ relationship
between teams, rather than an infinite (or zero) ratio of strength
parameters.

\subsection{Team-specific order effects}
\label{sectshfa}

In Section~\ref{secbtthree}, it is noted that team-specific home field
advantage parameters can be modelled by considering each team's home
performances independently of its road performances. The separation
properties of the data can therefore be determined exactly as in
Section~\ref{secestone} with $2t$ items, the home and road versions of
teams $1,2,\ldots,t$.

\subsection{Single order effect}\label{secsingleord}
If in (\ref{eqhomeprob}) we let $\beta_l=\log \pi_l$ and $\eta=\log
\gamma$, the Albert-Anderson discrimination condition for a game
between $i$ and $j$ with $i$ at home and winning (losing) is
$$ \beta_i-\beta_j+\eta \ge (\le) 0. $$
If we reparameterize, defining $\beta_{lH}=\beta_l+\eta/2$ and
$\beta_{lV}=\beta_l- \eta/2$, the condition becomes
\begin{equation}
\beta_{iH} - \beta_{jV} \ge (\le) 0.
\label{eqhomefirst}
\end{equation}
(H represents ``home team'' and V ``visiting team''.)
 This is the condition of
Section~\ref{sectshfa}, but since we have overparameterized, we also
have that 
\begin{equation}
\beta_{1H}-\beta_{1V}= \cdots = \beta_{tH}-\beta_{tV},
\label{eqhomeres}
\end{equation}
since each difference is equal to $\eta$. In inferring
$\trianglerighteq$ relationships between the teams, we therefore need
to take into account (\ref{eqhomeres}) as well as (\ref{eqhomefirst}).
Equation (\ref{eqhomeres}) has two implications:
\begin{itemize}
\item If there exists a set of teams $\{k_1,\ldots,k_r\}$, with $r\ge
2$, such that $\beta_{k_1H} \ge \beta_{k_2V}, \; \ldots \;
\beta_{k_{r-1}H} \ge \beta_{k_rV}, \; 
\beta_{k_rH} \ge \beta_{k_1V}$, 
then $\beta_{lH} \ge \beta_{lV}$ for all teams $l$, since $\eta \ge
0$. (The difference between the sum of the left-hand sides of the
inequalities and the sum of the right-hand sides is $r\eta$.)
The same applies with $H$ and $V$ interchanged.
\item If $\beta_{iH} \ge \beta_{jH}$ for any teams $i,j$, 
then $\beta_{iV} \ge \beta_{jV}$ for those teams
and conversely, since $\beta_i \ge \beta_j$.
\end{itemize}

The procedure of Section~\ref{secestone} is therefore no longer
sufficient. An iterative scheme will be necessary in its place, since
knowledge that $\eta \ge (\le) 0$ provides additional information
about whether $\beta_{iH} \ge (\le) \beta_{iV}$, and conversely. Note,
however, that if we once determine that $\eta \ge 0$, we do not need
to check again (since the implications for $\beta_{iH}$ and
$\beta_{iV}$ are fixed); the same applies separately when $0 \ge
\eta$.

To determine whether $\eta \ge 0$, we need to find a chain of
inequalities, as described above, beginning and ending at the same
team. Unless we already have that $\beta_{kH} \ge \beta_{kV}$ for some
$k$, this is most easily done by constructing a directed graph $G_H$
which has vertices representing the teams $\{1,2,\ldots,t\}$ and edges
from $i$ to $j$ whenever $\beta_{iH} \ge \beta_{jV}$ (that is, when
$iH \tre jV$). In $G_H$, no vertex $k$ is considered connected to
itself, since we only construct $G_H$ if $kH \tre kV$ does not hold
for any $k$. We then compute the transitive closure $\bar{G}_H$ of
$G_H$, and determine whether it contains any edges from a vertex $k$
to itself, in which case $\beta_{kH} \ge \beta_{kV}$ and hence $\eta
\ge 0$.

To determine whether $\eta \le 0$, the same technique can be applied
with $H$ and $V$ interchanged.

This leads to the following algorithm:
\begin{enumerate}
\item $lH \trianglerighteq lH$ and $lV \trianglerighteq lV$ for each
  team $l$ (for completeness).
\item For each game, let $i$ be the home team and $j$ the visiting
  team. If $i$ defeated $j$, $iH \trianglerighteq jV$; otherwise, $jV
  \trianglerighteq iH$.
\item \label{steptwo} \label{stepthree} Construct a directed graph $G$
  with the $2t$ vertices $1H,\ldots,tH,1V,\ldots,tV$ and edges
  connecting those vertices related by $\tre$. Compute the transitive
  closure $\bar{G}$ of $G$; any edges $(k,l)$ in $\bar{G}$ but not $G$
  imply an additional $\tre$ relationship between items $k$ and $l$.
\item \label{stepfive} If $kH \tre kV$ for any $k$, set $kH \tre kV$
  for all $k$, and go to Step \ref{stepsix}. Otherwise, construct a
  directed graph $G_H$ with vertices $1,2,\ldots,t$ and edges from $i$
  to $j$ whenever $iH \tre jV$.  Compute its transitive closure
  $\bar{G}_H$. If $\bar{G}_H$ contains an edge from some vertex $l$ to
  itself, set $kH \tre kV$ for all $k$.
\item \label{stepsix} If $kV \tre kH$ for any $k$, set $kV \tre kH$
  for all $k$, and go to Step \ref{stepfour}. Otherwise, construct a
  directed graph $G_V$ with vertices $1,2,\ldots,t$ and edges from $i$
  to $j$ whenever $iV \tre jH$.  Compute its transitive closure
  $\bar{G}_V$. If $\bar{G}_V$ contains an edge from some vertex $l$ to
  itself, set $kV \tre kH$ for all $k$.
\item \label{stepfour} For each $i$ and $j$, if $iH \tre jH$, then $iV
  \tre jV$, and conversely.
\item \label{stepeight} If steps \ref{steptwo} through \ref{stepfour}
  on the current cycle added any $\tre$ relationships, go back to
  step~\ref{steptwo}; otherwise, stop.
\end{enumerate}

It is computationally most convenient to maintain the directed graph
$G$ as an adjacency matrix, in which the $\tre$ relationships can be
found and updated directly.

The maximum number of $\tre$ relationships is $4t^2$, so the maximum
number of iterations of steps \ref{steptwo} through \ref{stepfour} is
bounded.  Often, in practice, only one or two iterations are
necessary.

To determine the probabilities whose maximum likelihood estimates are
1, 0 or arbitrary, we find which relation, $\gg, \cong, \gtrless$
applies between each pair of items, where a item is an entity such as
$iH$, a team with a venue attached.  (Note that even if $iH \gg jV$,
it is perfectly possible to have $jH \cong iV$, depending on the
magnitude of the home field advantage.)  Then, instead of using
(\ref{eqhomeprob}), we define the probability of $i$ defeating $j$
with $i$ at home as
$$ p_{ij} = \left\{
\begin{array}{rcl}
\displaystyle{\frac{\gamma\pi_i}{\gamma\pi_i+\pi_j}} 
& \mbox{if } & iH \cong jV\\
1 & \mbox{if } & iH \gg jV\\
0 & \mbox{if } & jV \gg iH\\
\mbox{ arbitrary} & \mbox{if } & iH \gtrless jV.
\end{array}
\right.
$$

Once again, using this modified probability formula in the maximum
likelihood equations effectively removes games whose outcomes have
maximum likelihood probabilities of 0 or 1 from the data set, and
produces a solution in which $\hat\gamma$ and all the $\hat\pi_i$
exist.

\section{Estimation in the presence of ties}\label{sectie}

\subsection{Reparameterization}\label{sectieone}

The discrimination conditions for models containing ties can be
written most conveniently if we define parameters which cast the two
versions of the model (single or team-specific tie parameters) in the
same notation, by defining $\pi_{i+}=\sqrt{\pi_i\nu}$ in the model
(\ref{eqtieprob}) or $\pi_{i+}=\sqrt{\pi_i\nu_i}$ in
(\ref{alltieprob}), so that the relative probabilities of a win, tie,
or loss for team $i$ in a game against team $j$ are $\pi_i$,
$\pi_{i+}\pi_{j+}$, and $\pi_j$, respectively. (Note that the model
with team-specific tie parameters could in fact be defined in this
way, with the $\pi_{i+}$ taking the place of the $\nu_i$ as
independent parameters.)  For convenience, we also define
$\pi_{i-}=\pi_i/\pi_{i+}$, which means that the relative probabilities
of a win or tie for team $i$ against team $j$ are $\pi_{i-}$ and
$\pi_{j+}$, respectively.

\citeN{albeande84} show that, if the response variable in a logistic
regression has $g$ categories, then each observation contributes $g-1$
conditions on their $\alpha$.  In the presence of ties, therefore,
each game gives two conditions on the $\beta_i$. If we write
$\beta_i=\log \pi_i$ and $\beta_{i+}=\log \pi_{i+}$ in either of the
models with ties, the conditions are
\begin{eqnarray*}
\beta_i - \beta_j  &\ge& 0 \\
\beta_i - \beta_{i+} - \beta_{j+} &\ge& 0 
\end{eqnarray*}
if $i$ defeats $j$, and
\begin{eqnarray*}
\beta_{i+}+\beta_{j+}-\beta_i &\ge& 0 \\
\beta_{i+}+\beta_{j+}-\beta_j &\ge& 0 
\end{eqnarray*}
if the result is a tie.

If we define an additional parameter
$\beta_{i-}=\beta_i-\beta_{i+}=\log\pi_{i-}$, the conditions take on
the simpler form
\begin{equation}
\begin{array}{rcl}
\beta_i  &\ge& \beta_j \\
\beta_{i-} &\ge& \beta_{j+}
\end{array}
\label{eqconda}
\end{equation}
if $i$ defeats $j$, and
\begin{equation}
\begin{array}{rcl}
\beta_{j+}  &\ge& \beta_{i-} \\
\beta_{i+}  &\ge& \beta_{j-}
\end{array}
\label{eqcondd}
\end{equation}
if the result is a tie.

To determine which results have maximum-likelihood probabilities which
are 0 or arbitrary, we need to find the exhaustive set of inequalities
which apply among the $t$ parameters $\beta_1,\ldots,\beta_t$ and also
among the $2t$ parameters
$\beta_{1+},\ldots,\beta_{t+},\beta_{1-},\ldots,\beta_{t-}$ as a
result of the conditions (\ref{eqconda}) and (\ref{eqcondd}), along
with auxiliary conditions relating $\beta_i$, $\beta_{i+}$, and
$\beta_{i-}$.  These will be different depending on whether the model
contains a single tie parameter or one for each team.

\subsection{Single tie parameter}
\label{secsingletie}

If there is only one tie parameter, the definitions of the parameters
imply $\beta_{i\pm}=(\beta_i\pm\log\nu)/2$.  This means that the
discrimination conditions can be re-cast without reference to
$\beta_i$ by replacing (\ref{eqconda}) with the conditions
$$
\beta_{i\pm}\ge\beta_{j\pm}
$$
by adding or subtracting $\log\nu$ from both sides.
The auxiliary relationship
$$
\beta_{1+}-\beta_{1-} = \cdots = \beta_{t+}-\beta_{t-}
$$
is identical to the condition (\ref{eqhomeres}) given in
Section~\ref{secsingleord}, so the procedure determining the complete
set of inequalities applying to this model is the same as that in the
model with a single order effect, even though it is being applied in a
different context. Specifically:
\begin{enumerate}
\item $l+ \tre l+$ and $l- \tre l-$ for each
team $l$ (for completeness).
\item For each game in which $i$ defeated $j$, $i+\tre j+$, $i-\tre
  j-$, and $i-\tre j+$.  For each game in which $i$ tied $j$, $i+\tre
  j-$ and $j+\tre i-$.
\item[\ref{stepthree}-\ref{stepeight}.] Carry out the corresponding
  steps from the single-order-parameter algorithm of
  Section~\ref{secsingleord}, with $H$ replaced by $+$ and $V$ by $-$.
\end{enumerate}

Once the full set of $\tre$ relationships has been determined, they
can be used to determine the $\gg$, $\cong$, or $\gtrless$
relationships between the $2t$ items $1+,\ldots,t+,1-,\ldots,t-$.
Since $i+$ is guaranteed by step \ref{stepfour} of the algorithm to
have the same relationship with $j+$ as $i-$ does with $j-$ (and to
allow an analogy with the following section), it is convenient to
refer to that relationship as holding between $i$ and $j$.  For
example, if $i+\gg j+$ (which means $i-\gg j-$) we say $i\gg j$.

{From} these relationships, we can conclude that the maximum
likelihood estimates of certain probabilities are 0, 1 or arbitrary.
However, with three possible outcomes, we may have one outcome whose
estimated probability is zero, while the remaining two probabilities
are both estimated to be nonzero. Tables~\ref{tabfourone} and
\ref{tabfourtwo} show, respectively, the maximum-likelihood
probabilities of wins and ties taking the relationships into account.
The maximum likelihood estimates
$\hat\pi_i,\hat\pi_j,\hat\pi_{i+},\hat\pi_{j+}$ are guaranteed to
exist.

\begin{table}
\begin{center}
  \begin{tabular}{r|cccc|}
    $p_{ij1}$ & $i\gg j$ & $i\cong j$ & $j\gg i$ & $i\gtrless j$ \\
    \hline
    $i-\gg j+$ & 1 & $\frac{\pi_i}{\pi_i+\pi_j}$ & 0 & arbitrary \\
    $i-\cong j+$ & $\frac{\pi_i}{\pi_i+\pi_{i+}\pi_{j+}}$ 
    & $\frac{\pi_i}{\pi_i+\pi_j+\pi_{i+}\pi_{j+}}$ & 0 
    & arbitrary \\
    $j+\gg i-$ & 0 & 0 & 0 & 0 \\
    $i-\gtrless j+$ & arbitrary & arbitrary & 0 & arbitrary \\
    \hline
  \end{tabular}
\caption{\label{tabfourone}Revised definition of $p_{ij1}$,
  the probability of team $i$ defeating team $j$}
\end{center}
\end{table}

\begin{table}
\begin{center}
  \begin{tabular}{r|cccc|}
    $p_{ij0}$ & $i+\gg j-$ & $i+\cong j-$ & $j-\gg i+$ 
    & $i+\gtrless j-$ \\
    \hline
    $j+\gg i-$ & 1 & $\frac{\pi_{i+}\pi_{j+}}{\pi_j+\pi_{i+}\pi_{j+}}$
    & 0 & arbitrary \\
    $j+\cong i-$ & $\frac{\pi_{i+}\pi_{j+}}{\pi_i+\pi_{i+}\pi_{j+}}$ 
    & $\frac{\pi_{i+}\pi_{j+}}{\pi_i+\pi_j+\pi_{i+}\pi_{j+}}$ 
    & 0 & arbitrary \\
    $i-\gg j+$ & 0 & 0 & 0 & 0 \\
    $j+\gtrless i-$ & arbitrary & arbitrary & 0 & arbitrary \\
    \hline
  \end{tabular}
\end{center}
\caption{\label{tabfourtwo}Revised definition of $p_{ij0}$,
the probability of a tie between teams $i$ and $j$}
\end{table}

In the Tables, note that the combination of relationships which gives
an arbitrary probability for any result is only possible if the two
teams have not actually played. Note also that in some cases (for
example $i\cong j$, $i+\gtrless j-$, $i-\gtrless j+$) the ratio of the
probabilities of two outcomes can be well-defined and finite, with the
ratio of either to the probability of the third outcome unconstrained
by the data. That is, the conditional probabilities of the first and
second outcomes, given that the third outcome does not happen, lie
strictly between 0 and 1. However, in this case we still consider all
three probabilities to be arbitrary.

If the maximum likelihood estimates are zero for two of the
probabilities $p_{ij0},p_{ij1},p_{ij2}$, the third probability is
estimated to be 1, and any such games can simply be removed from the
data set before estimating the remaining probabilities.  However, in
games where exactly one outcome is constrained to have zero
probability, the modified probabilities from Table~\ref{tabfourone}
and \ref{tabfourtwo} must be used to leave only those outcomes out of
likelihood maximization.  Since the likelihood equations here equate
observed and expected points for each team (where points are defined
as in Section~\ref{secbtfour}), as well as observed and expected
overall ties, this means that we must restrict our calculations to
those outcomes whose probabilities are known to lie strictly between 0
and 1.

If we make the obvious generalization that $y_{ijk}$ is the number of
games between $i$ and $j$ with outcome $k$, and take $p_{ijk}$ from
Tables~\ref{tabfourone} and \ref{tabfourtwo}, the likelihood equations
become, for each $\pi_i$,
$$
2 \sum_{\substack{j\\ 0<p_{ij1}<1}} y_{ij1} 
+ \sum_{\substack{j\\ 0<p_{ij0}<1}} y_{ij0} 
= 2 \sum_{\substack{j\\ 0<p_{ij1}<1}} n_{ij} p_{ij1} 
+ \sum_{\substack{j\\ 0<p_{ij0}<1}} n_{ij}
p_{ij0}
,
$$
and for $\nu$,
$$
\sum_{\substack{i<j\\ 0<p_{ij1}<1}} y_{ij0} 
= \sum_{\substack{i<j\\ 0<p_{ij1}<1}} n_{ij} p_{ij0}
.
$$
These equations yield a method analogous to Ford's for finding the
maximum likelihood estimates.  In the first set of equations, a factor
$\sqrt{\pi_i}$ can be removed from the right-hand side and isolated,
and in the last equation the same can be done with $\nu$. Then, in
each case, a fixed-point iteration can be carried out.

\subsection{Team-specific tie parameters}\label{sectietwo}

In the case of a team-specific tie parameter model, the auxiliary
equations describing the overparameterization that exists due to the
introduction of the extra $t$ parameters
$\beta_{1-},\ldots,\beta_{t-}$ are
$$
\beta_{i}=\beta_{i+}+\beta_{i-}
.
$$
These do not lend themselves to a direct simplification of
(\ref{eqconda}), but allow us to infer some inequalities among
parameters by adding and subtracting others.  For example, if
$\beta_{i+}\ge\beta_{j-}$ and $\beta_{i-}\ge\beta_{j+}$, we can add
the two inequalities to infer $\beta_i\ge\beta_j$.  This produces the
following algorithm for determining the $\tre$ relationships amongst
the teams $1,\ldots,t$ as well as the items
$1+,\ldots,t+,1-,\ldots,t-$:
\begin{enumerate}
\item $l \tre l$, $l+ \tre l+$, and $l-
  \tre l-$ for each team $l$ (for completeness).
\item For each game in which $i$ defeated $j$, $i\tre j$ and $i-\tre
  j+$.  For each game in which $i$ tied $j$, $i+\tre j-$ and $j+\tre
  i-$.
\item \label{loopstart} Construct a directed graph $G$ with vertices
  $1,2,\ldots,t$ and with edges from $i$ to $j$ if $i \tre j$. Compute
  the transitive closure $\bar{G}$ of $G$. Any edges $(k,l)$ in
  $\bar{G}$ but not $G$ imply additional relationships $k \tre l$.
  Repeat with the directed graph $G+$ whose edges are the $2t$ items
  $1+,\dots,t+,1-,\ldots,t-$.
\item \label{loopend} For each pair of teams $i$ and $j$:
  \begin{itemize}
  \item If $i\tre j$ and $j+\tre i+$, then $i-\tre j-$;
  \item If $i\tre j$ and $j+\tre i-$, then $i+\tre j-$;
  \item If $i\tre j$ and $j-\tre i+$, then $i-\tre j+$;
  \item If $i\tre j$ and $j-\tre i-$, then $i+\tre j+$;
  \item If $i+\tre j+$ and $i-\tre j-$, then $i\tre j$;
  \item If $i+\tre j-$ and $i-\tre j+$, then $i\tre j$.
  \end{itemize}
\item If steps \ref{loopstart} or \ref{loopend} on the current cycle
  added any $\tre$ relationships, go back to step~\ref{loopstart};
  otherwise, stop.
\end{enumerate}
As in Section~\ref{secsingleord}, it is computationally convenient to
maintain $G$ and $G+$ as adjacency matrices, for ease of determination
and updating of the $\tre$ relationships.

Once we have obtained the full set of $\tre$ relationships among the
$t$ teams and the $2t$ team-and-sign items, we use them to find
the $\cong$, $\gg$, and $\gtrless$ relationships among teams and
between the ``plussed'' teams $1+,\ldots,t+$ and the ``minussed''
teams $1-,\ldots,t-$.  These can then be used to identify the
outcomes whose maximum-likelihood probabilities vanish
and to deduce correctly the remaining estimated probabilities
$p_{ijk}$ according to Tables~\ref{tabfourone} and \ref{tabfourtwo}.
Again, this procedure permits a maximum
likelihood solution for the remaining probabilities 
with parameter estimates that exist.  This time,
the likelihood equations are, for each $\pi_i$,
$$
\sum_{\substack{j\\ 0<p_{ij1}<1}} y_{ij1}
 = \sum_{\substack{j\\ 0<p_{ij1}<1}} n_{ij} p_{ij1}
,
$$
and for each $\nu_i$,
$$
\sum_{\substack{j\\ 0<p_{ij1}<1}}
 y_{ij0} = \sum_{\substack{j\\ 0<p_{ij1}<1}} n_{ij} p_{ij0}
.
$$
Again, a Ford-like iterative process can be constructed. In the
first set of equations, a $\pi_i$ can be factored out of the
right-hand side and isolated, and in the second set, the same can be
done with either $\sqrt{\nu_i}$ or $\pi_{i+}$.

\section{Presentation of results}\label{secpres}

\subsection{Introduction}
Our procedures yield maximum likelihood estimates for all the
probabilities $p_{ij}$ (or $p_{ijk}$ in the presence of ties), but
these probabilities may be determined from the existence of, say, a
$\gg$ relationship rather than from finite parameters that represent
the strengths of the teams. It is, nonetheless, desirable to have a
``strength'' measure for each team regardless of whether the parameter
estimates all exist; for example, we may wish to produce a rank
ordering of the teams.

\subsection{Simple Bradley-Terry model}\label{secpresone}

Here, the maximum-likelihood solution is completely specified by a
partition of the teams into equivalence classes according to the
$\cong$ relation defined in Section~\ref{secestone}, the relations
($\gg$ or $\gtrless$) between classes, and the estimated $\pi_i$ for
each team, defined up to an overall multiplicative factor \emph{within
  each class}.  This does not permit easy comparison
of teams in different classes.

To provide a ranking of the teams even in the presence of $\gtrless$
relationships, as well as a single number with a concrete
interpretation, we can define the round-robin winning percentage, or
``RRWP''.  This is the proportion $R_i$ of games a team would be
expected to win if they played each other team an equal number of
times (as in a ``balanced schedule'' or ``round-robin tournament''),
and is given by the average of that team's estimated winning
probability against all of the other teams:
$$
R_i=\frac{1}{t-1}\sum_{j\ne i} p_{ij}.
$$
Arbitrary $p_{ij}$ are (for definiteness) set equal to
$\frac{1}{2}$ in this calculation.  This allows the RRWP to be defined
given any set of results.  This definition gives a sensible ordering
of the teams, since if $p_{ij}>p_{ji}$ for a pair of teams $i$ and
$j$, then (because $p_{ik}\ge p_{jk}$ for any other team $k$),
$R_i>R_j$.  If two teams are related by $\gtrless$, and thus the
maximum likelihood probability of one defeating the other is
arbitrary, the RRWP can be used to provide a somewhat arbitrary
ordering of them.  We do not have enough information to compare such
teams directly in the Bradley-Terry model, but the RRWP in effect
pulls to the middle of the rankings teams with many arbitrary
$p_{ij}$, those teams being ``unknown quantities''.

In the case where the teams have actually played a balanced schedule,
each team's RRWP will be equal to its actual winning percentage.  This
is a manifestation of the familiar result that the Bradley-Terry model
reproduces ordering by won-lost record in the case of a balanced
schedule (``complete block design'' in the language of
\citeN{bradterr52}).  However, this definition extends to the case
where some teams are related by $\gg$, in which case there are maximum
likelihood estimates of parameters in the traditional Bradley-Terry
model that do not exist.

\subsection{Order effects}\label{secprestwo}

In the presence of order effects, either single or team-specific, the
definition of RRWP must be generalized to include the home field
advantage.  The quantity $p_{ij}$ is now the predicted probability for
$i$ to win a game at home against $j$, and because of the home field
advantage $p_{ji}\ne 1-p_{ij}$ in general.  Note that the ``diagonal''
probability $p_{ii}$ is a measure of the strength of the home
advantage; with a single order effect, it is simply
$\gamma/(1+\gamma)$, while with team-specific order effects, $p_{ii}$
measures the difference between team $i$'s performances as a home team
and as a visiting team.

A logical generalization of the RRWP to a model with order effects is
the percentage of games a team would be expected to win if they played
each other team once at home and once on the road:
\begin{equation}\label{eqrrwptwo}
R_i=\frac{1}{2(t-1)}\sum_{j\ne i} (p_{ij} + 1 - p_{ji});
\end{equation}
Once again, arbitrary $p_{ij}$ are replaced in the sum by 1/2, and if
the teams have actually played this kind of round-robin, each team's
RRWP will be equal to its actual winning percentage.  In a model with
a single order effect where the home advantage parameter exists
($kH\cong kV$ for all $k$), $R_i>R_j$ if $i\gg j$, or if $i\cong j$
and $\pi_i>\pi_j$.

\subsection{Models with ties}
\label{secpresthree}

When ties are allowed in the model, the concept of RRWP must be
generalized.  Our procedures now give estimated probabilities
$p_{ijk}$ of each result (win, loss or tie) in games between each pair
of teams.  We can use the $p_{ijk}$ to obtain the predicted number of
wins, losses, and ties for each team in a round-robin schedule.  The
predicted round-robin percentage $R_{ik}$ of outcome $k$ for team $i$
is
\begin{equation}
  R_{ik}=\frac{1}{(t-1)}\sum_{j\ne i} p_{ijk}
  .
\end{equation}
As in Sections~\ref{secbtfour} and \ref{secbtfive}, $k=1$ denotes a
win for team $i$, $k=2$ a win for team $j$, and $k=0$ a tie.  All
arbitrary probabilities for a particular hypothetical (unplayed) game
are assigned equal weight, so that if all three probabilities (win,
loss, or tie) against a particular opponent are arbitrary, they are
all assumed to be $1/3$ in the sum, while if one probability $p_{ijk}$
vanishes, the other two are assumed to be $1/2$.  One could of course
define a more sophisticated method of making this assignment, taking
into consideration known ratios of unknown probabilities, or the
overall proportion of ties, but we only present the simplest method
here.

These round-robin proportions can be used to define a RRWP, according
to the traditional definition of a tie as a half a win and a half a
loss, by
\begin{equation}
  R_{i}=R_{i1}+\frac{1}{2}R_{i0}
  .
\end{equation}
In the single tie parameter model, if $p_{ij1}>p_{ij2}$, then
$R_i>R_j$.  This means that if no teams are related by $\gtrless$ the
ranking of the teams by strength and by RRWP will be identical.  This
is \emph{not} true in the model with team-specific tie parameters, as
illustrated in Section~\ref{secexrank}.

The round-robin proportions predicted by the team-specific tie
parameter model will always agree with the actual fractions of games
won, lost, and tied by each team after a balanced schedule has been
completed.  For the model with a single tie parameter, it is only the
combination wins-plus-half-losses (or equivalent) which has this
property. Thus, with team-specific tie parameters, there is nothing
special about counting a tie as half a win; indeed, given any system
awarding $c_k$ points for outcome $k$, we can define a ``round-robin
points per game'', RRPPG, that is equivalent to RRWP if the point
system is a linear transformation of $c_1=1,c_2=0,c_0=\frac{1}{2}$:
$$
T_i=\sum_{k=0}^2 c_k R_{ik}
.
$$
After a round-robin tournament with this point system, the points
obtained per game by each team will be equal to $T_i$. This implies
that the $T_i$ may sensibly be used to rank the teams after any
collection of games. 

The world standard for soccer, for example, sets $c_1=3,c_2=0,c_0=1$
-- that is, three points for a win and one for a tie.  A model with a
single tie parameter (or the idea, mentioned in
Section~\ref{secbtfour}, of counting ties literally as half a win and
half a loss) is no longer reasonable here, since the teams could be
ranked differently by $T_i$ and by points even at the end of a
balanced schedule (because the observed and expected proportions of
ties for each team may not match). Instead, it is necessary to fit
team-specific tie parameters and {\em then} rank the teams by $T_i$,
since then observed points per game will equal $T_i$ after a balanced
schedule.

\section{Examples}\label{secex}

\subsection{Example 1}

Suppose there are four teams, $a,b,c,d$. Teams $a$ and $b$ play twice,
winning one game each, and $c$ and $d$ do the same. Team $a$ plays $c$
once, with $a$ winning. Even though each team has at least one win and
at least one defeat, it is seen that there is still a separation
between the groups $\{a,b\}$ and $\{c,d\}$, because the only
inter-group comparison is $a$'s win over $c$.

Applying the procedure of Section~\ref{secestone} gives the relations
$a\ge b,\; b \ge a,\; c \ge d,\; d \ge c,\; a \ge c$. Hence $a \tre
b,\; a \tre c,\; b \tre a,\; c \tre d,\; d \tre c$ directly and $a
\tre d, \; b \tre c, \; b \tre d$ indirectly, from computing the
transitive closure. These are all the $\tre$ relationships, and so the
teams are related as follows:
$$
\begin{array}{cccc}
a \\
\cong & b\\
\gg & \gg & c\\
\gg & \gg & \cong & d
\end{array}
$$

In other words, $a$ and $b$ as a group are ``infinitely stronger''
than $c$ and $d$ as a group, and the maximum likelihood probability of
a team in $\{a,b\}$ defeating a team in $\{c,d\}$ is 1.

The maximum likelihood estimates for the probabilities of $a$
defeating $b$ and for $c$ defeating $d$ are both (unsurprisingly) 0.5,
and so the RRWP's for $a$, $b$, $c$ and $d$ are 0.833, 0.833, 0.167
and 0.167 respectively.

\subsection{Example 2}

Suppose four teams play the following games. The first-named team is
at home in each case, and we will want to fit a home field advantage. 
\begin{itemize}
\item $a$ vs $b$: $a$ wins
\item $b$ vs $a$: $b$ wins
\item $c$ vs $a$: $a$ wins
\item $c$ vs $d$: $c$ wins
\item $d$ vs $c$: $d$ wins
\end{itemize}
(This is in fact the same set of games as in Example 1.)

The procedure of Section~\ref{secsingleord} is concisely expressed in
the table below. A number in a particular row and column denotes a
$\tre$ relationship between the items of that row and column; the
number itself denotes the number of the step at which that $\tre$
relationship was added. An asterisk denotes that the $\tre$
relationship was added on the second iteration. A third iteration
failed to add any more $\tre$ relationships.

\begin{center}
\begin{tabular}{r|cccc|cccc|}
& $aH$ & $bH$ & $cH$ & $dH$ & $aV$ & $bV$ & $cV$ & $dV$ \\
\hline
$aH$ & 1 & & 6* & 6 & 4 & 2 & 3* & 3*\\
$bH$ & & 1 & 3 & & 2 & 4 & 3* & 3\\
$cH$ & & & 1 & & & & 4 & 2\\
$dH$ & & & & 1 & & & 2 & 4\\
\hline
$aV$ & & & 2 & & 1 & & 3* & 3 \\
$bV$ & & & & & & 1 & 6 &\\
$cV$ & & & & & & & 1 & \\
$dV$ & & & & & & & & 1\\
\hline
\end{tabular}
\end{center}

The relationships between each $H$-team and each $V$-team are
therefore (the $\ll$ relation between $CH$ and $AV$ of course means
$AV\gg CH$):

\begin{center}
\begin{tabular}{r|cccc|r}
& $aV$ & $bV$ & $cV$ & $dV$ & RRWP\\
\hline
$aH$ & $\gg$ & $\gg$ & $\gg$ & $\gg$ & 0.750\\
$bH$ & $\gg$ & $\gg$ & $\gg$ & $\gg$ & 0.667\\
$cH$ & $\ll$ & $\gtrless$ & $\gg$ & $\gg$ & 0.250\\
$dH$ & $\gtrless$ & $\gtrless$ & $\gg$ & $\gg$ & 0.333\\
\hline
\end{tabular}
\end{center}

The RRWP's apply to each of teams $a,b,c,d$ as a whole (averaged over
the items representing their home and road performances, as in
(\ref{eqrrwptwo})).

The estimated probabilities of the events that occurred are all 1, and
hence the maximum value of the likelihood is 1. This result contrasts
with Example~1, which was based on the same data, but without a home
field effect. The home field advantage in this Example is
overwhelming, to the extent that the winner of a game between $a$ and
$b$ is predicted ``with certainty'' to be the home team, but the
difference in strength between $a$ and $c$ is so large that it even
overcomes this home field advantage.  Finally, since $aV \gg cH$ but
$bV \gtrless cH$, the RRWP for team $a$ is greater than that for team
$b$, as shown.

\subsection{Example 3}
\label{secexrank}

Suppose now that three teams $a$, $b$, $c$ play a round-robin in which
each team meets the others four times, as follows:

\begin{center}
\begin{tabular}{rl}
$a$ vs.\ $b$: & $b$ wins once, 3 ties\\
$a$ vs.\ $c$: & $a$ wins 4 times\\
$b$ vs.\ $c$: & $c$ wins 2 times, 2 ties.
\end{tabular}
\end{center}

The won-lost-tied records of these three teams are as follows:

\begin{center}
\begin{tabular}{rl}
$a$: & won 4, lost 1, tied 3\\
$b$: & won 1, lost 2, tied 5\\
$c$: & won 2, lost 4, tied 2
\end{tabular}
\end{center}

With a single tie parameter, the algorithm of
Section~\ref{secsingletie} shows that all items $a+,b+,\ldots,c-$ are
related $\cong$: the likelihood is maximized by probabilities strictly
between 0 and 1. The estimated probabilities are:


\begin{center}
\begin{tabular}{rr|rrr}
$i$ & $j$ & $p_{ij1}$ & $p_{ij2}$ & $p_{ij0}$\\
\hline
$a$ & $b$ & 0.464 & 0.126 & 0.410 \\
$a$ & $c$ & 0.513 & 0.101 & 0.385 \\
$b$ & $c$ & 0.316 & 0.229 & 0.455 
\end{tabular}
\end{center}
and the round-robin wins,
losses, and ties per game, and RRWP, are:

\begin{center}
\begin{tabular}{r|rrr|r}
Team, $i$
& $R_{i1}$ & $R_{i2}$ & $R_{i0}$  & $R_i$ \\
\hline
$a$ &  0.489 & 0.114 & 0.398 & 0.6875\\
$b$ &  0.221 & 0.346 & 0.432 & 0.4375\\
$c$ &  0.164 & 0.414 & 0.420 & 0.3750
\end{tabular}
\end{center}

Team $b$ is estimated to be stronger than team $c$, since
$p_{bc1}>p_{bc2}$, which is not surprising considering that $b$ has 7
points in the round-robin (2 points for a win and 1 for a tie), while
$c$ has only 6.  Note that the RRWP agrees with the actual winning
percentage for each of the three teams, but the round-robin wins,
losses, and ties per game for each team predicted by the model do not
agree with the actual numbers.

If we apply the algorithm of Section~\ref{sectietwo} to fit
team-specific tie parameters, we obtain the $\tre$ relationships shown
in the two tables below.  A number in a table indicates, as
previously, that the teams or items in question are related by $\tre$,
with the number itself indicating the step at which the $\tre$ was
found.

\begin{center}
\begin{tabular}{ccc}
\begin{tabular}{r|rrr|}
& $a$ & $b$ & $c$\\
\hline
$a$ & 1 & 3 & 2\\
$b$ & 2 & 1 & 3\\
$c$ & 3 & 2 & 1\\
\hline
\end{tabular}
&&
\begin{tabular}{r|ccc|ccc|}
& $a+$ & $b+$ & $c+$ & $a-$ & $b-$ & $c-$\\
\hline
$a+$ & 1 & & 4 && 2 &\\
$b+$ & 3&1&3&2&3&2\\
$c+$ & 3&&1&&2&\\
\hline
$a-$& 3&4&2&1&3&4\\
$b-$&2&&4&&1&\\
$c-$&3&2&3&3&3&1\\
\hline
\end{tabular}
\end{tabular}
\end{center}

This shows that $a- \gg c+$ and $c- \gg a+$, with all other relevant
relationships being $\cong$. (Only relationships of the form $i$ vs.\ 
$j$, $i+$ vs.\ $j-$ with $i\ne j$ concern games that can be played.)
{From} the second Table of Section~\ref{secsingletie}, the estimated
probability of a tie between $a$ and $c$ is zero. The entire set of
estimated probabilities is:


\begin{center}
\begin{tabular}{rr|rrr}
$i$ & $j$ & $p_{ij1}$ & $p_{ij2}$ & $p_{ij0}$ \\
\hline
$a$ & $b$ & 0.210 & 0.040 & 0.750\\
$a$ & $c$ & 0.790 & 0.210 & 0 \\
$b$ & $c$ & 0.210 & 0.290 & 0.500
\end{tabular}
\end{center}
and the RRWP values are:

\begin{center}
\begin{tabular}{r|rrr|r}
Team, $i$
& $R_{i1}$ & $R_{i2}$ & $R_{i0}$ & $R_i$\\
\hline
$a$ & 0.500 & 0.125 & 0.375 & 0.6875\\
$b$ & 0.125 & 0.250 & 0.625 & 0.4375\\
$c$ & 0.250 & 0.500 & 0.250 & 0.3750
\end{tabular}
\end{center}

Note now that team $c$ is estimated to be stronger than $b$, since its
probability of winning against $b$ is higher than its probability of
losing. In other words, under the team-specific tie parameter model,
the ranking of teams by strength no longer necessarily agrees with the
ranking under a point system of 2 per win and 1 per tie. The
explanation here is that $b$ has almost no chance of defeating $a$,
but gains a substantial number of points against $a$ from ties --
enough, in fact, to overcome its deficit of points from games against
$c$.

Note also that in this model the predicted proportions of round-robin
wins, losses, and ties for each team agree with the observed
proportions.

\subsection{Example 4}
\label{secexwac}


\renewcommand{\baselinestretch}{1.0}
\begin{table}[htbp]
  \begin{center}
    \begin{tabular}{l|l|l||rr|rr|rr}
      &&& \multicolumn{2}{c|}{Oct 18}&
      \multicolumn{2}{c|}{Oct 25}&
      \multicolumn{2}{c}{Final}\\
      Team & Abbr.\ & Div.\ & W-L & RRWP & W-L & RRWP& W-L & RRWP\\
      \hline
      Air Force & AFA & Mtn & 3-1 & .667 & 4-1 & .635 & 7-1 & .849\\
      San Diego State 
      & SDSU & Pac & 3-0 & .600 & 4-0 & .967 & 7-1 & .821\\
      Brigham Young
      & BYU & Pac & 2-1 & .233 & 3-1 & .488 & 7-1 & .813\\
      Wyoming
      & Wyo.\ & Mtn & 3-0 & .967 & 4-0 & .967 & 6-2 & .764\\
      Rice
      & Rice & Mtn & 2-1 & .667 & 2-2 & .450 & 5-3 & .629\\
      Colorado State
      & CSU & Mtn & 3-1 & .667 & 4-1 & .658 & 5-3 & .629\\
      Utah & Utah & Pac & 2-1 & .600 & 2-2 & .867 & 5-3 & .613\\
      Fresno State
      & FSU & Pac & 1-2 & .300 & 2-2 & .609 & 5-3 & .576\\
      Texas Christian
      & TCU & Mtn & 2-1 & .667 & 2-2 & .610 & 4-4 & .560\\
      Southern Methodist
      & SMU & Mtn & 2-2 & .667 & 3-2 & .529 & 4-4 & .499\\
      Texas-El Paso
      & UTEP & Pac &2-1 & .667 & 2-2 & .441 & 3-5 & .337 \\
      San Jose State 
      & SJSU & Pac & 2-1 & .667 & 2-2 & .380 & 3-5 & .337\\
      Tulsa & Tulsa & Mtn & 0-3 & .200 & 0-4 & .100 & 2-6 & .372\\
      New Mexico
      & UNM & Pac & 0-4 & .200 & 1-4 & .133 & 1-7 & .100\\
      Nevada-Las Vegas
      & UNLV & Mtn & 0-4 & .133 & 0-5 & .100 & 0-8 & .067\\
      Hawaii
      & Haw.\ & Pac & 0-4 & .100 & 0-5 & .067 & 0-8 & .033\\
      \hline
    \end{tabular}
    \caption{Summary of results from WAC data
      of Section~\ref{secexwac}}
    \label{tab:kenwac}
  \end{center}
\end{table}
\renewcommand{\baselinestretch}{1.8}

\begin{figure}
  \begin{center}
    \begin{picture}(130,240)(-70,-240)
      \put(0,-10){\oval(30,20)}
      \put(-13,-13){Wyo.}
      \put(-10,-20){\line(0,-1){20}}
      \put(-10,-20){\vector(0,-1){12}} 
      \put(10,-20){\line(1,-2){25}}
      \put(10,-20){\vector(1,-2){16}}  
      \put(-35,-50){\oval(60,20)}
      \put(-50,-53){Others}
      \put(-10,-60){\line(3,-4){30}} 
      \put(-10,-60){\vector(3,-4){17}} 
      \put(-15,-60){\line(0,-1){90}} 
      \put(-15,-60){\vector(0,-1){47}} 
      \put(-55,-60){\line(0,-1){90}} 
      \put(-55,-60){\vector(0,-1){47}} 
      \put(-35,-80){\oval(30,20)}
      \put(-48,-83){SDSU}
      \put(-45,-90){\line(0,-1){60}} 
      \put(-45,-90){\vector(0,-1){32}} 
      \put(-35,-90){\line(0,-1){110}} 
      \put(-35,-90){\vector(0,-1){57}} 
      \put(-25,-90){\line(0,-1){60}} 
      \put(-25,-90){\vector(0,-1){32}} 
      \put(45,-80){\oval(30,20)}
      \put(32,-83){Utah}
      \put(37.5,-90){\line(0,-1){10}}       
      \put(37.5,-90){\vector(0,-1){7}} 
      \put(30,-110){\oval(30,20)}
      \put(17,-113){FSU}
      \put(30,-120){\line(0,-1){10}} 
      \put(30,-120){\vector(0,-1){7}} 
      \put(30,-140){\oval(30,20)}
      \put(17,-143){BYU}
      \put(25,-150){\line(-1,-1){50}} 
      \put(25,-150){\vector(-1,-1){27}} 
      \put(37.5,-150){\line(0,-1){20}} 
      \put(37.5,-150){\vector(0,-1){12}} 
      \put(-55,-160){\oval(30,20)}
      \put(-68,-163){Tulsa}
      \put(-15,-160){\oval(30,20)}
      \put(-28,-163){UNM}
      \put(45,-180){\oval(30,20)}
      \put(31,-183){UNLV}
      \put(-35,-210){\oval(30,20)}
      \put(-48,-213){Haw.}
\put(-20,-240){Figure \protect\ref{fig:kenrels}a}
\end{picture}
\hfill
\begin{picture}(70,190)(-70,-170)
      \put(-55,10){\oval(30,20)}
      \put(-68,7){Wyo.}
      \put(-45,0){\line(0,-1){10}}
      \put(-45,0){\vector(0,-1){7}}  
      \put(-15,10){\oval(30,20)}
      \put(-28,7){SDSU}
      \put(-25,0){\line(0,-1){10}}
      \put(-25,0){\vector(0,-1){7}}  
      \put(-35,-20){\oval(30,20)}
      \put(-48,-23){Utah}
      \put(-35,-30){\line(0,-1){10}}       
      \put(-35,-30){\vector(0,-1){7}} 
      \put(-35,-50){\oval(60,20)}
      \put(-50,-53){Others}
      \put(-35,-60){\line(0,-1){10}} 
      \put(-35,-60){\vector(0,-1){5}} 
      \put(-15,-60){\line(0,-1){40}} 
      \put(-15,-60){\vector(0,-1){22}} 
      \put(-55,-60){\line(0,-1){40}} 
      \put(-55,-60){\vector(0,-1){22}} 
      \put(-35,-80){\oval(30,20)}
      \put(-48,-83){UNM}
      \put(-35,-90){\line(0,-1){40}} 
      \put(-35,-90){\vector(0,-1){22}} 
      \put(-55,-110){\oval(30,20)}
      \put(-68,-113){Tulsa}
      \put(-15,-110){\oval(30,20)}
      \put(-28,-113){UNLV}
      \put(-35,-140){\oval(30,20)}
      \put(-48,-143){Haw.}
\put(-60,-170){Figure \protect\ref{fig:kenrels}b}
    \end{picture}
\hfill
\begin{picture}(60,130)(-65,-170)
      \put(-35,-50){\oval(60,20)}
      \put(-50,-53){Others}
      \put(-50,-60){\line(0,-1){10}} 
      \put(-50,-60){\vector(0,-1){5}} 
      \put(-20,-60){\line(0,-1){40}} 
      \put(-20,-60){\vector(0,-1){22}} 
      \put(-50,-80){\oval(30,20)}
      \put(-63,-83){UNM}
      \put(-50,-90){\line(0,-1){40}} 
      \put(-50,-90){\vector(0,-1){22}} 
      \put(-20,-110){\oval(30,20)}
      \put(-33,-113){UNLV}
      \put(-50,-140){\oval(30,20)}
      \put(-63,-143){Haw.}
\put(-60,-170){Figure \protect\ref{fig:kenrels}c}
\end{picture}
    \caption{Diagrams illustrating the relationships among teams at
      three points during the 1998 WAC Football season.}
  \label{fig:kenrels}
  \end{center}
\end{figure}

In Fall 1998, the Western Athletic Conference (WAC), an American
collegiate sports association, had 16 members. For (American) football
competition, these were divided geographically into two divisions.
Each team played all the other teams in its own division, but faced
only one opponent from the other division. The teams, their
abbreviated names and their division membership are shown in the first
three columns of Table~\ref{tab:kenwac}. ``Mtn'' denotes Mountain
division, ``Pac'' Pacific. Game results were obtained from
\citeN{usatoday98}.

We focus on three subsets of the game results: the games played before
October 18, those played before October 25, and all the games played
in the regular season. There were no ties, and we have chosen not to
model any home field advantage.  Model (\ref{btorig}) was fitted to
each subset of the data, and RRWP calculated for each team as in
Section~\ref{secpresone}.  Table~\ref{tab:kenwac} summarizes the wins,
losses and RRWP for each team and subset.

The RRWP values do not show which relationships, $\cong$, $\gg$, $\ll$
or $\gtrless$, hold between each pair of teams. These are depicted in
Figure~\ref{fig:kenrels} for each of the three subsets of the data.
The relation $\cong$, as described in Section~\ref{secestone}, divides
the teams into equivalence classes; in this example, all the classes
except one contain only one team (which is used to label the class),
while the remaining class, labelled ``Others'', contains all the
remaining teams. The vertical position of the classes on the page is
determined by RRWP; an arrow connecting two classes denotes a $\gg$
relationship between them. (For simplicity, arrows implied by
transitivity are not shown; thus, for example, Wyo.\ $\gg$ UNM in all
three diagrams.)

Most teams played only three games before October 18. As a result,
there are numerous $\gg$ and $\gtrless$ relationships among the teams,
and the relationship diagram, shown in Figure~\ref{fig:kenrels}a, is
complicated. However, the diagram explains two apparent anomalies in
the RRWP values from Table~\ref{tab:kenwac}: SDSU's three wins came
against teams that lost all their games, explaining the low RRWP of
.600, while BYU's low RRWP of .233 comes from defeating two winless
teams (Hawaii and UNLV) and losing to FSU, which was FSU's only win.

\renewcommand{\baselinestretch}{1.0}
\begin{table}[htbp]
  \begin{center}
    \begin{tabular}{ll|cc}
      &&\multicolumn{2}{c}{Relationship}\\
      Winner & Loser & Before & After\\
      \hline
      Air Force & Tulsa & $\gg$ & $\gg$\\
      Wyoming & Rice & $\gg$ & $\gg$ \\
      CSU & TCU & $\cong$ & $\cong$ \\
      SMU & UNLV & $\gg$ & $\gg$ \\
      SDSU & Utah & $\gtrless$ & $\gg$ \\
      BYU & SJSU & $\ll$ & $\cong$ \\
      FSU & UTEP & $\ll$ & $\cong$ \\
      UNM & Hawaii & $\gtrless$ & $\gg$
    \end{tabular}
    \caption{Game results Oct 18--Oct 25 and effect on relationships}
    \label{tab:weektwo}
  \end{center}
\end{table}
\renewcommand{\baselinestretch}{1.8}

Each team played once between October 18 and October 25, with the
results shown in Table~\ref{tab:weektwo}. Results where the winning
team was already related to the losing team by $\gg$ or $\cong$ do not
affect the relationships, but a relationship of $\gtrless$ changes to
$\gg$ and one of $\ll$ changes to $\cong$, pending the results of
other games. Thus BYU and FSU join SJSU and UTEP in the ``Others''
class, while now SDSU $\gg$ Utah and UNM $\gg$ Hawaii. This simplifies
the relationship diagram considerably, as shown in
Figure~\ref{fig:kenrels}b. Note that Utah now has the third-highest
RRWP because the team's two losses came against the top two teams.

At the end of the season, the relationship diagram was as shown in
Figure~\ref{fig:kenrels}c. Thirteen of the sixteen teams are now
related $\cong$, so that estimated win probabilities for games between
these teams are strictly between 0 and 1, and can be estimated with
standard techniques. UNM, UNLV and Hawaii are still related $\ll$ to
the ``Others'' class, since none of these teams have any wins against
teams in the ``Others'' class. Because UNM $\gg$ Hawaii but UNLV
$\gtrless$ UNM and Hawaii (that is, UNM's only win of the season came
against Hawaii, while UNLV played neither team), the three teams'
RRWP's rank as shown in the last column of Table~\ref{tab:kenwac}.

The final RRWPs for the teams in the ``Others'' class generally agree
with the teams' won-lost records. Differences in RRWP between teams
with the same numbers of wins arise for two reasons: the Mountain
division is stronger overall than the Pacific, and some teams faced a
stronger opponent in their inter-divisional game than others. For
example, Air Force (Mountain) came out ahead of SDSU and BYU
(Pacific), while SDSU was ahead of BYU because their inter-divisional
opponents were Tulsa and UNLV respectively. On the other hand, Rice
and CSU had identical 5-3 records and had interdivisional opponents
with identical 3-5 records, so both pairs of teams ended with
identical RRWP's.

\section{Discussion}

In the Bradley-Terry model, we have shown that the separation
properties of the data, and therefore the determination of the
probabilities estimated to be 0 or 1 by maximum likelihood, can be
found by considering the transitive closure of certain graphs. This
makes for straightforward computation, and avoids the need for the
linear programming routines that are required for logistic regression
in general.  We also suggest that the Round-Robin Winning Percentage
(RRWP) or Round-Robin Points Per Game (RRPPG) provide a sensible way
of summarizing the results, because they are defined regardless of the
separation properties of the data.

One practical application of the RRWP (or RRPPG) is to provide
``mid-season'' rankings for a league of teams which will eventually
play a balanced schedule, but have not yet completed it.  At the end
of the season, the each team's RRWP will be equal to the actual
winning percentage used by the league to order the final standings,
but when only part of the season has been played, the RRWP should be a
more accurate measure of a team's performance to date, as it considers
the strength of the completed portion of their schedule.  Another
application would be to a league which, because of considerations such
as travel or length of season, does not play a balanced schedule.
Traditionally, the teams in such leagues are ranked according to
winning percentage, which benefits teams who happen to play weaker
opponents more often.

We also note that in models with team-specific order effects (such as
(\ref{allhomeprob}) or team-specific tie parameters (such as
(\ref{alltieprob})), there is more than one parameter per team, so
that ranking the teams by a single strength parameter is either
impossible (in the former case) or not necessarily desirable (in the
latter).  RRWP and RRPPG provide a uniform means for ranking which may
be preferable depending on the application.

\section*{Aknowledgments}

We would like to thank G.~Hatfield and R.~Stagat for useful
discussions of applications of the Bradley-Terry model.  JTW
acknowledges support by the Swiss Nationalfonds, and by the Tomalla
Foundation, Z\"{u}rich.

\bibliography{myref}

\end{document}